\documentclass[12pt,reqno]{amsart}
\usepackage{amscd,amssymb,amsmath,multicol,float,scalefnt,cancel,array}
\restylefloat{table}
\newtheorem{thm}[equation]{Theorem}
\numberwithin{equation}{section}

\begin{document}
\setcounter{MaxMatrixCols}{26}
\raggedbottom \voffset=-.7truein \hoffset=0truein \vsize=8truein
\hsize=6truein \textheight=8truein \textwidth=6truein
\baselineskip=18truept

\def\mapright#1{\ \smash{\mathop{\longrightarrow}\limits^{#1}}\ }
\def\mapleft#1{\smash{\mathop{\longleftarrow}\limits^{#1}}}
\def\mapup#1{\Big\uparrow\rlap{$\vcenter {\hbox {$#1$}}$}}
\def\mapdown#1{\Big\downarrow\rlap{$\vcenter {\hbox {$\ssize{#1}$}}$}}
\def\mapne#1{\nearrow\rlap{$\vcenter {\hbox {$#1$}}$}}
\def\mapse#1{\searrow\rlap{$\vcenter {\hbox {$\ssize{#1}$}}$}}
\def\mapr#1{\smash{\mathop{\rightarrow}\limits^{#1}}}
\def\ss{\smallskip}
\def\vp{v_1^{-1}\pi}
\def\at{{\widetilde\alpha}}
\def\sm{\wedge}
\def\la{\langle}
\def\ra{\rangle}
\def\on{\operatorname}
\def\ol#1{\overline{#1}{}}
\def\spin{\on{Spin}}
\def\lbar{\ell}
\def\qed{\quad\rule{8pt}{8pt}\bigskip}
\def\ssize{\scriptstyle}
\def\a{\alpha}
\def\bz{{\Bbb Z}}
\def\im{\on{im}}
\def\ct{\widetilde{C}}
\def\ext{\on{Ext}}
\def\sq{\on{Sq}}
\def\eps{\epsilon}
\def\ar#1{\stackrel {#1}{\rightarrow}}
\def\br{{\bold R}}
\def\bC{{\bold C}}
\def\bA{{\bold A}}
\def\bB{{\bold B}}
\def\bD{{\bold D}}
\def\bh{{\bold H}}
\def\bQ{{\bold Q}}
\def\bP{{\bold P}}
\def\bx{{\bold x}}
\def\bo{{\bold{bo}}}
\def\si{\sigma}
\def\Vbar{{\overline V}}
\def\dbar{{\overline d}}
\def\wbar{{\overline w}}
\def\Sum{\sum}
\def\tfrac{\textstyle\frac}
\def\tb{\textstyle\binom}
\def\Si{\Sigma}
\def\w{\wedge}
\def\equ{\begin{equation}}
\def\b{\beta}
\def\G{\Gamma}
\def\g{\gamma}
\def\k{\kappa}
\def\psit{\widetilde{\Psi}}
\def\tht{\widetilde{\Theta}}
\def\psiu{{\underline{\Psi}}}
\def\thu{{\underline{\Theta}}}
\def\aee{A_{\text{ee}}}
\def\aeo{A_{\text{eo}}}
\def\aoo{A_{\text{oo}}}
\def\aoe{A_{\text{oe}}}
\def\vbar{{\overline v}}
\def\endeq{\end{equation}}
\def\sn{S^{2n+1}}
\def\zp{\bold Z_p}
\def\A{{\cal A}}
\def\P{{\mathcal P}}
\def\cQ{{\mathcal Q}}
\def\cj{{\cal J}}
\def\zt{{\bold Z}_2}
\def\bs{{\bold s}}
\def\bof{{\bold f}}
\def\bq{{\bold Q}}
\def\be{{\bold e}}
\def\Hom{\on{Hom}}
\def\ker{\on{ker}}
\def\coker{\on{coker}}
\def\da{\downarrow}
\def\colim{\operatornamewithlimits{colim}}
\def\zphat{\bz_2^\wedge}
\def\io{\iota}
\def\Om{\Omega}
\def\Prod{\prod}
\def\e{{\cal E}}
\def\exp{\on{exp}}
\def\abar{{\overline a}}
\def\xbar{{\overline x}}
\def\ybar{{\overline y}}
\def\zbar{{\overline z}}
\def\Rbar{{\overline R}}
\def\nbar{{\overline n}}
\def\cbar{{\overline c}}
\def\qbar{{\overline q}}
\def\bbar{{\overline b}}
\def\et{{\widetilde E}}
\def\ni{\noindent}
\def\coef{\on{coef}}
\def\den{\on{den}}
\def\lcm{\on{l.c.m.}}
\def\vi{v_1^{-1}}
\def\ot{\otimes}
\def\psibar{{\overline\psi}}
\def\mhat{{\hat m}}
\def\exc{\on{exc}}
\def\ms{\medskip}
\def\ehat{{\hat e}}
\def\etao{{\eta_{\text{od}}}}
\def\etae{{\eta_{\text{ev}}}}
\def\dirlim{\operatornamewithlimits{dirlim}}
\def\gt{\widetilde{L}}
\def\lt{\widetilde{\lambda}}
\def\st{\widetilde{s}}
\def\ft{\widetilde{f}}
\def\sgd{\on{sgd}}
\def\lfl{\lfloor}
\def\rfl{\rfloor}
\def\ord{\on{ord}}
\def\gd{{\on{gd}}}
\def\rk{{{\on{rk}}_2}}
\def\nbar{{\overline{n}}}
\def\lg{{\on{lg}}}
\def\cB{\mathcal{B}}
\def\cS{\mathcal{S}}
\def\cP{\mathcal{P}}
\def\N{{\Bbb N}}
\def\Z{{\Bbb Z}}
\def\Q{{\Bbb Q}}
\def\R{{\Bbb R}}
\def\C{{\Bbb C}}
\def\l{\left}
\def\r{\right}
\def\mo{\on{mod}}
\def\xt{\times}
\def\notimm{\not\subseteq}
\def\Remark{\noindent{\it  Remark}}

\def\*#1{\mathbf{#1}}
\def\0{$\*0$}
\def\1{$\*1$}
\def\22{$(\*2,\*2)$}
\def\33{$(\*3,\*3)$}
\def\ss{\smallskip}
\def\ssum{\sum\limits}
\def\dsum{\displaystyle\sum}
\def\la{\langle}
\def\ra{\rangle}
\def\on{\operatorname}
\def\o{\on{o}}
\def\U{\on{U}}
\def\lg{\on{lg}}
\def\a{\alpha}
\def\bz{{\Bbb Z}}
\def\eps{\varepsilon}
\def\bc{{\bold C}}
\def\bN{{\bold N}}
\def\nut{\widetilde{\nu}}
\def\tfrac{\textstyle\frac}
\def\b{\beta}
\def\G{\Gamma}
\def\g{\gamma}
\def\zt{{\Bbb Z}_2}
\def\zth{{\bold Z}_2^\wedge}
\def\bs{{\bold s}}
\def\bx{{\bold x}}
\def\bof{{\bold f}}
\def\bq{{\bold Q}}
\def\be{{\bold e}}
\def\lline{\rule{.6in}{.6pt}}
\def\xb{{\overline x}}
\def\xbar{{\overline x}}
\def\ybar{{\overline y}}
\def\zbar{{\overline z}}
\def\ebar{{\overline \be}}
\def\nbar{{\overline n}}
\def\rbar{{\overline r}}
\def\Mbar{{\overline M}}
\def\et{{\widetilde e}}
\def\ni{\noindent}
\def\ms{\medskip}
\def\ehat{{\hat e}}
\def\what{{\widehat w}}
\def\Yhat{{\widehat Y}}
\def\nbar{{\overline{n}}}
\def\minp{\min\nolimits'}
\def\mul{\on{mul}}
\def\N{{\Bbb N}}
\def\Z{{\Bbb Z}}
\def\Q{{\Bbb Q}}
\def\R{{\Bbb R}}
\def\C{{\Bbb C}}
\def\notint{\cancel\cap}
\def\cS{\mathcal S}
\def\cR{\mathcal R}
\def\el{\ell}
\def\TC{\on{TC}}
\def\dstyle{\displaystyle}
\def\ds{\dstyle}
\def\Wbar{\wbar}
\def\zcl{\on{zcl}}
\def\Vb#1{{\overline{V_{#1}}}}

\def\Remark{\noindent{\it  Remark}}
\title[Approach to TC of Klein bottle]
{An approach to the topological complexity of the Klein bottle}
\author{Donald M. Davis}
\address{Department of Mathematics, Lehigh University\\Bethlehem, PA 18015, USA}
\email{dmd1@lehigh.edu}
\date{February 1, 2017}


\maketitle
\begin{abstract} Recently, Cohen and Vandembroucq proved that the reduced topological complexity of the Klein bottle is 4. Simultaneously and independently we announced a proof of the same result. Mistakes were found in our argument, which was quite different than theirs. After correcting these, we found that our description of the obstruction class agreed with theirs.
Our approach to showing that this obstruction is nonzero failed to do so, while theirs did not fail. Here we discuss our approach, which deals more directly with the simplicial structure of the Klein bottle.
 \end{abstract}
\section{Introduction}
The reduced topological complexity, $\TC(X)$, of a topological space $X$, as introduced in \cite{F}, is roughly one less than the minimal number of rules required to tell how to move between any two points of $X$. A ``rule'' must be such that the choice of path varies continuously with the endpoints.
An outstanding problem has been to determine $\TC(K)$, where $K$ is the Klein bottle.

In a recent preprint (\cite{CV}), Cohen and Vandembroucq  proved that $\TC(K)=4$. Simultaneously and independently we announced a proof of the same result. Mistakes were found in our argument, which was quite different than theirs. After correcting these, we found that our description of the obstruction class agreed with theirs.
Our approach to showing that this obstruction is nonzero failed to do so, while theirs did not fail. Here we discuss our approach, which deals more directly with the simplicial structure of the Klein bottle.

The Klein bottle is homeomorphic to the space of all configurations of various physical systems. For example, it is homeomorphic to the space of all planar 5-gons with side lengths 1, 1, 2, 2, and 3 (\cite[Table B]{pps}). Such a polygon can be considered as linked robot arms. Knowing that $\TC(K)=4$ implies that five rules are required to program these arms to move from any configuration to any other.

Our main tool is a result of Costa and Farber (\cite{CF}), which we state later as Theorem \ref{CFthm}, which describes a single obstruction in $H^{2n}(X\times X;G)$, where $G$ is a certain local coefficient system, for an $n$-dimensional cell complex $X$ to satisfy $\TC(X)< 2n$. We present an approach to proving that this class is nonzero for the Klein bottle $K$, which would imply that $\TC(K)\ge4$, while $\TC(K)\le4$ for dimensional reasons.(\cite[Cor 4.15]{F})

\section{The $\Delta$-complex for $K\times K$}\label{Deltasec}
A $\Delta$-complex, as described in \cite{Ha}, is essentially a quotient of a simplicial complex, with certain simplices identified. As discussed there, this notion is equivalent to that of semi-simplicial complex introduced in \cite{EZ}. It is important that vertices be numbered prior to identifications, and that simplices be described by writing vertices in increasing order.

The $\Delta$-complex that we will use for $K$ is given below. It has one vertex $v$, three edges, $(0,2)=(4,5)$, $(1,2)=(3,4)$, and $(0,1)=(3,5)$, and two 2-cells, $(0,1,2)$ and $(3,4,5)$.

\bigskip
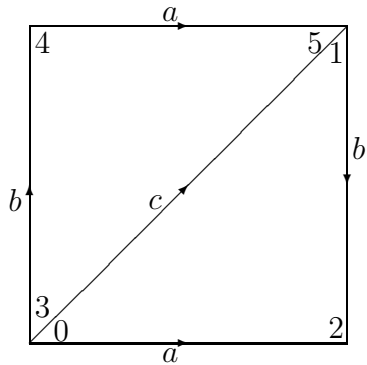
\begin{center}
\begin{picture}(120,125)
\put(0,0){\vector(1,0){60}}
\put(60,0){\line(1,0){60}}
\put(120,120){\vector(0,-1){60}}
\put(120,60){\line(0,-1){60}}
\put(0,0){\vector(0,1){60}}
\put(0,60){\line(0,1){60}}
\put(0,120){\vector(1,0){60}}
\put(60,120){\line(1,0){60}}
\put(0,0){\vector(1,1){60}}
\put(60,60){\line(1,1){60}}
\put(50,-7){$a$}
\put(9,1){$0$}
\put(50,122){$a$}
\put(122,70){$b$}
\put(-8,50){$b$}
\put(45,51){$c$}
\put(2,10){$3$}
\put(2,110){$4$}
\put(113,2){$2$}
\put(113,106){$1$}
\put(105,110){$5$}
\end{picture}
\end{center}

\bigskip
If $K$ and $L$ are simplicial complexes with an ordering of the vertices of each, then the simplices of the simplicial complex $K\times L$ are all $\la(v_{i_0},w_{j_0}),\ldots,(v_{i_k},w_{j_k})\ra$ such that $i_0\le \cdots\le i_k$ and $j_0\le\cdots\le j_k$ and $\{v_{i_0},\ldots,v_{i_k}\}$ and $\{w_{j_0},\ldots,w_{j_k}\}$ are simplices of $K$ and $L$, respectively. Note that we may have $v_{i_t}=v_{i_{t+1}}$ or $w_{j_t}=w_{j_{t+1}}$, but not both (for the same $t$). Now, if $K$ and $L$ are $\Delta$-complexes, i.e., they have some simplices identified, then $K\times L$ has
$$\la(v_{i_0},w_{j_0}),\ldots,(v_{i_k},w_{j_k})\ra\sim\la(v_{i_0'},w_{j_0'}),\ldots,(v_{i_k'},w_{j_k'})\ra$$
iff $\{v_{i_0},\ldots,v_{i_k}\}\sim\{v_{i_0'},\ldots,v_{i_k'}\}$ and $\{w_{j_0}\ldots w_{j_k}\}\sim\{w_{j_0'}\ldots w_{j_k'}\}$, and the positions of the repetitions in $(v_{i_0},\ldots,v_{i_k})$ and $(v_{i_0'},\ldots,v_{i_k'})$ are the same, and so are those of $(w_{j_0},\ldots,w_{j_k})$ and $(w_{j_0'},\ldots,w_{j_k'})$. This description is equivalent to the one near the end of \cite{EZ2}, called $K\Delta L$ there. There is also a discussion in \cite[pp.277-278]{Ha}.

Following this description, we now list the simplices of $K\times K$, where $K$ is the above $\Delta$-complex. We write $v$ for the unique vertex when it is being producted with a simplex, but otherwise we list all vertices by their number. We omit commas in ordered pairs; e.g., $24$ denotes the vertex of $K\times K$ which is vertex 2 in the first factor and vertex 4 in the second factor. There are 1, 15, 50, 60, and 24 distinct simplices of dimensions 0, 1, 2, 3, and 4, respectively. This is good since $1-15+50-60+24=0$, the Euler characteristic of $K\times K$.
We number the simplices in each dimension, which will be useful later.

The only 0-simplex is $(vv)$.

\medskip
1-simplices:

1: $(0v,1v)=(3v,5v)$.

2: $(0v,2v)=(4v,5v)$.

3: $(1v,2v)=(3v,4v)$.

4: $(v0,v1)=(v3,v5)$.

5: $(v0,v2)=(v4,v5)$.

6: $(v1,v2)=(v3,v4)$.

7: $(00,11)=(30,51)=(33,55)=(03,15)$.

8: $(00,12)=(30,52)=(34,55)=(04,15)$.

9: $(01,12)=(31,52)=(33,54)=(03,14)$.

10: $(00,21)=(40,51)=(43,55)=(03,25)$.

11: $(00,22)=(40,52)=(44,55)=(04,25)$.

12: $(01,22)=(41,52)=(43,54)=(03,24)$.

13: $(10,21)=(30,41)=(33,45)=(13,25)$.

14: $(10,22)=(30,42)=(34,45)=(14,25)$.

15: $(11,22)=(31,42)=(33,44)=(13,24)$.

\medskip
2-simplices

1: $(0v,1v,2v)$.

2: $(3v,4v,5v)$.

3: $(v0,v1,v2)$.

4: $(v3,v4,v5)$.

5: $(00,10,11)=(30,50,51)=(03,13,15)=(33,53,55)$.

6: $(00,10,12)=(30,50,52)=(04,14,15)=(34,54,55)$.

7: $(01,11,12)=(31,51,52)=(03,13,14)=(33,53,54)$.

8: $(00,01,11)=(30,31,51)=(03,05,15)=(33,35,55)$.

9: $(00,02,12)=(30,32,52)=(04,05,15)=(34,35,55)$.

10: $(01,02,12)=(31,32,52)=(03,04,14)=(33,34,54)$.

11: $(00,20,21)=(40,50,51)=(03,23,25)=(43,53,55)$.

12: $(00,20,22)=(40,50,52)=(04,24,25)=(44,54,55)$.

13: $(01,21,22)=(41,51,52)=(03,23,24)=(43,53,54)$.

14: $(00,01,21)=(40,41,51)=(03,05,25)=(43,45,55)$.

15: $(00,02,22)=(40,42,52)=(04,05,25)=(44,45,55)$.

16: $(01,02,22)=(41,42,52)=(03,04,24)=(43,44,54)$.

17: $(10,20,21)=(30,40,41)=(13,23,25)=(33,43,45)$.

18: $(10,20,22)=(30,40,42)=(14,24,25)=(34,44,45)$.

19: $(11,21,22)=(31,41,42)=(13,23,24)=(33,43,44)$.

20: $(10,11,21)=(30,31,41)=(13,15,25)=(33,35,45)$.

21: $(10,12,22)=(30,32,42)=(14,15,25)=(34,35,45)$.

22: $(11,12,22)=(31,32,42)=(13,14,24)=(33,34,44)$.

23: $(00,11,22)$.

24: $(30,41,52)$.

25: $(03,14,25)$.

26: $(33,44,55)$.

27: $(00,01,22)=(40,41,52)$.

28: $(30,31,52)=(00,01,12)$.

29: $(03,04,25)=(43,44,55)$.

30: $(33,34,55)=(03,04,15)$.

31: $(00,10,22)=(04,14,25)$.

32: $(30,40,52)=(34,44,55)$.

33: $(03,13,25)=(00,10,21)$.

34: $(33,43,55)=(30,40,51)$.

35: $(00,21,22)=(40,51,52)$.

36: $(03,24,25)=(43,54,55)$.

37: $(30,51,52)=(00,11,12)$.

38: $(33,54,55)=(03,14,15)$.

39: $(00,12,22)=(04,15,25)$.

40: $(30,42,52)=(34,45,55)$.

41: $(03,15,25)=(00,11,21)$.

42: $(33,45,55)=(30,41,51)$.

43: $(01,12,22)=(03,14,24)$.

44: $(31,42,52)=(33,44,54)$.

45: $(03,13,24)=(01,11,22)$.

46: $(33,43,54)=(31,41,52)$.

47: $(10,21,22)=(30,41,42)$.

48: $(13,24,25)=(33,44,45)$.

49: $(10,11,22)=(30,31,42)$.

50: $(13,14,25)=(33,34,45)$.

\medskip
3-simplices

1: $(00,10,21,22)$

2: $(00,11,21,22)$.

3: $(00,10,11,22)$.

4: $(00,11,12,22)$.

5: $(00,01,11,22)$.

6: $(00,01,12,22)$.

7: $(30,40,51,52)$.

8: $(30,41,51,52)$.

9: $(30,40,41,52)$.

10: $(30,41,42,52)$.

11: $(30,31,41,52)$.

12: $(30,31,42,52)$.

13: $(03,13,24,25)$.

14: $(03,14,24,25)$.

15: $(03,13,14,25)$.

16: $(03,14,15,25)$.

17: $(03,04,14,25)$.

18: $(03,04,15,25)$.

19: $(33,43,54,55)$.

20: $(33,44,54,55)$.

21: $(33,43,44,55)$.

22: $(33,44,45,55)$.

23: $(33,34,44,55)$.

24: $(33,34,45,55)$.

25: $(00,10,20,21)=(03,13,23,25)$.

26: $(00,10,11,21)=(03,13,15,25)$.

27: $(00,01,11,21)=(03,05,15,25)$.

28: $(00,10,20,22)=(04,14,24,25)$.

29: $(00,10,12,22)=(04,14,15,25)$.

30: $(00,02,12,22)=(04,05,15,25)$.

31: $(01,11,21,22)=(03,13,23,24)$.

32: $(01,11,12,22)=(03,13,14,24)$.

33: $(01,02,12,22)=(03,04,14,24)$.

34: $(00,01,02,12)=(30,31,32,52)$.

35: $(00,01,11,12)=(30,31,51,52)$.

36: $(00,10,11,12)=(30,50,51,52)$.

37: $(00,01,02,22)=(40,41,42,52)$.

38: $(00,01,21,22)=(40,41,51,52)$.

39: $(00,20,21,22)=(40,50,51,52)$.

40: $(10,11,12,22)=(30,31,32,42)$.

41: $(10,11,21,22)=(30,31,41,42)$.

42: $(10,20,21,22)=(30,40,41,42)$.

43: $(30,40,50,51)=(33,43,53,55)$.

44: $(30,40,41,51)=(33,43,45,55)$.

45: $(30,31,41,51)=(33,35,45,55)$.

46: $(30,40,50,52)=(34,44,54,55)$.

47: $(30,40,42,52)=(34,44,45,55)$.

48: $(30,32,42,52)=(34,35,45,55)$.

49: $(31,41,51,52)=(33,43,53,54)$.

50: $(31,41,42,52)=(33,43,44,54)$.

51: $(31,32,42,52)=(33,34,44,54)$.

52: $(03,04,05,15)=(33,34,35,55)$.

53: $(03,04,14,15)=(33,34,54,55)$.

54: $(03,13,14,15)=(33,53,54,55)$.

55: $(03,04,05,25)=(43,44,45,55)$.

56: $(03,04,24,25)=(43,44,54,55)$.

57: $(03,23,24,25)=(43,53,54,55)$.

58: $(13,14,15,25)=(33,34,35,45)$.

59: $(13,14,24,25)=(33,34,44,45)$.

60: $(13,23,24,25)=(33,43,44,45)$.

\medskip
4-simplices

1: $(00,10,20,21,22)$.

2: $(00,10,11,21,22)$.

3: $(00,10,11,12,22)$.

4: $(00,01,11,21,22)$.

5: $(00,01,11,12,22)$.

6: $(00,01,02,12,22)$.

7: $(30,40,50,51,52)$.

8: $(30,40,41,51,52)$.

9: $(30,40,41,42,52)$.

10: $(30,31,41,51,52)$.

11: $(30,31,41,42,52)$.

12: $(30,31,32,42,52)$.

13: $(03,13,23,24,25)$.

14: $(03,13,14,24,25)$.

15: $(03,13,14,15,25)$.

16: $(03,04,14,24,25)$.

17: $(03,04,14,15,25)$.

18: $(03,04,05,15,25)$.

19: $(33,43,53,54,55)$.

20: $(33,43,44,54,55)$.

21: $(33,43,44,45,55)$.

22: $(33,34,44,54,55)$.

23: $(33,34,44,45,55)$.

24: $(33,34,35,45,55)$.

\section{$H^4(K\times K)$ with local coefficients}\label{H4sec}
We will need to show that a certain class is nonzero in $H^4(K\times K;G)$ with coefficients in a certain local coefficient system $G$. In this section, we describe the relations
in $H^4(K\times K;G)$ for an arbitrary free abelian local coefficient system $G$.

For a $\Delta$-complex $X$ with a single vertex $x_0$, such as the one just described for $K\times K$, a local coefficient system $G$ is an abelian group $G$ together with an action of  $\pi_1(X;x_0)$ on $G$, giving $G$ the structure of $\Z[\pi_1(X;x_0)]$-module. If $C_k$ denotes the free abelian group generated by the $k$-cells of $X$, homomorphisms
$$\delta_{k-1}:\Hom(C_{k-1},G)\to\Hom(C_k,G)$$
are defined by
\begin{equation}\label{delta}\delta_{k-1}(\phi)(\la v_{i_0},\ldots,v_{i_k}\ra)=\rho_{i_0,i_1}\cdot\phi(\la v_{i_1},\ldots,v_{i_k}\ra)+\sum_{i=1}^k(-1)^i\phi(\la v_{i_0},\ldots,\widehat{ v_i},\ldots,v_{i_k}\ra),\end{equation}
where $\widehat{v_i}$ denotes omission of that vertex, and $\rho_{i_0,i_1}$ is the element of $\pi_1(X;x_0)$ corresponding to the edge from $v_{i_0}$ to $v_{i_1}$.
Then
$$H^k(X;G)=\ker(\delta_k)/\im(\delta_{k-1}).$$
This description is given in \cite[p.501]{EZ}.

We have $$\pi_1(K;v)=\la a,b,c\ra/(c=ab^{-1}=ba)$$
and $$\pi_1(K\times K;(v,v))=\la a,b,c,a',b',c'\ra/(c=ab^{-1}=ba,\,c'=a'b'^{-1}=b'a'),$$
where the primes correspond to the second factor. We will prove the following key result, in which  $\g_1,\ldots, \g_{24}$ denote the generators corresponding to the 4-cells of $K\times K$
listed at the end of the previous section, and $G$ is any free abelian local coefficient system.
\begin{thm}\label{H4thm}  Let
$$\eps_j=\begin{cases}-1&j\equiv2\mod 3\\
1&j\equiv0,1\mod3,\end{cases}$$
and let $\psi\in\Hom(C_4,G)$.
Then $[\psi]=0\in H^4(K\times K;G)$ if and only if $\ds\sum_{j=1}^{24}\eps_j\psi(\g_j)$ is 0 in the quotient of $G$ modulo the action of $b-1$, $b'-1$, $c+1$, and $c'+1$.
\end{thm}
{\it Proof}.
We find the image of $\delta:\Hom(C_3,G)\to \Hom(C_4,G)$ by row-reducing, using only integer operations, the $60$-by-$24$ matrix $M$ whose entries $m_{i,j}\in\Z[\pi_1(K\times K)]$ satisfy
\begin{equation}\label{ij}\delta(\phi)(\g_j)=\sum_{i=1}^{60}m_{i,j}\phi(\b_i),\quad j=1,\ldots,24.\end{equation}
Here $\b_i$ denotes the generator of $C_3$ corresponding to the $i$th 3-cell.

We list the matrix $M$ below. Each row has two nonzero entries, as each 3-cell is a face of two 4-cells, while each column has five nonzero entries, as each 4-cell is bounded by five 3-cells. For example, $m_{42,1}=c$ because $(10,20,21,22)$ is obtained from $(00,10,20,21,22)$ by omission of the initial vertex, and $\rho_{00,10}=c$.

The row reduction of this matrix  can be done by hand in less than 30 minutes.
We  find that the row-reduced form  has 27 nonzero rows, with its only nonzero elements in rows $i=1,\ldots,23$ being 1 in position $(i,i)$ and $-\eps_i$ in position $(i,24)$, while in rows 24 through 27 the only nonzero element is in column 24, and equals $b-1$, $b'-1$, $c+1$, and $c'+1$, respectively.

A row $(r_1,\ldots,r_{24})$, with $r_j=\psi(\g_j)$ for $\psi\in\Hom(C_4,G)$, is equivalent, modulo the first 23 rows just described, to a row with 0's in the first 23 positions, and $\sum_{j=1}^{24}\eps_jr_j$ in the
24th column. The last four rows of the reduced matrix yield the claim of Theorem \ref{H4thm}.

{\scalefont{.85}{
$$\renewcommand\arraystretch{.75}
\arraycolsep=1.5pt
\begin{array}{r|cccccccccccccccccccccccc}
&1&2&3&4&5&6&7&8&9&10&11&12&13&14&15&16&17&18&19&20&21&22&23&24\\
\hline
1&1&1&0&0&0&0&0&0&0&0&0&0&0&0&0&0&0&0&0&0&0&0&0&0\\
2&0&-1&0&-1&0&0&0&0&0&0&0&0&0&0&0&0&0&0&0&0&0&0&0&0\\
3&0&-1&-1&0&0&0&0&0&0&0&0&0&0&0&0&0&0&0&0&0&0&0&0&0\\
4&0&0&-1&0&-1&0&0&0&0&0&0&0&0&0&0&0&0&0&0&0&0&0&0&0\\
5&0&0&0&-1&-1&0&0&0&0&0&0&0&0&0&0&0&0&0&0&0&0&0&0&0\\
6&0&0&0&0&1&1&0&0&0&0&0&0&0&0&0&0&0&0&0&0&0&0&0&0\\
7&0&0&0&0&0&0&1&1&0&0&0&0&0&0&0&0&0&0&0&0&0&0&0&0\\
8&0&0&0&0&0&0&0&-1&0&-1&0&0&0&0&0&0&0&0&0&0&0&0&0&0\\
9&0&0&0&0&0&0&0&-1&-1&0&0&0&0&0&0&0&0&0&0&0&0&0&0&0\\
10&0&0&0&0&0&0&0&0&-1&0&-1&0&0&0&0&0&0&0&0&0&0&0&0&0\\
11&0&0&0&0&0&0&0&0&0&-1&-1&0&0&0&0&0&0&0&0&0&0&0&0&0\\
12&0&0&0&0&0&0&0&0&0&0&1&1&0&0&0&0&0&0&0&0&0&0&0&0\\
13&0&0&0&0&0&0&0&0&0&0&0&0&1&1&0&0&0&0&0&0&0&0&0&0\\
14&0&0&0&0&0&0&0&0&0&0&0&0&0&-1&0&-1&0&0&0&0&0&0&0&0\\
15&0&0&0&0&0&0&0&0&0&0&0&0&0&-1&-1&0&0&0&0&0&0&0&0&0\\
16&0&0&0&0&0&0&0&0&0&0&0&0&0&0&-1&0&-1&0&0&0&0&0&0&0\\
17&0&0&0&0&0&0&0&0&0&0&0&0&0&0&0&-1&-1&0&0&0&0&0&0&0\\
18&0&0&0&0&0&0&0&0&0&0&0&0&0&0&0&0&1&1&0&0&0&0&0&0\\
19&0&0&0&0&0&0&0&0&0&0&0&0&0&0&0&0&0&0&1&1&0&0&0&0\\
20&0&0&0&0&0&0&0&0&0&0&0&0&0&0&0&0&0&0&0&-1&0&-1&0&0\\
21&0&0&0&0&0&0&0&0&0&0&0&0&0&0&0&0&0&0&0&-1&-1&0&0&0\\
22&0&0&0&0&0&0&0&0&0&0&0&0&0&0&0&0&0&0&0&0&-1&0&-1&0\\
23&0&0&0&0&0&0&0&0&0&0&0&0&0&0&0&0&0&0&0&0&0&-1&-1&0\\
24&0&0&0&0&0&0&0&0&0&0&0&0&0&0&0&0&0&0&0&0&0&0&1&1\\
25&1&0&0&0&0&0&0&0&0&0&0&0&-1&0&0&0&0&0&0&0&0&0&0&0\\
26&0&1&0&0&0&0&0&0&0&0&0&0&0&0&1&0&0&0&0&0&0&0&0&0\\
27&0&0&0&1&0&0&0&0&0&0&0&0&0&0&0&0&0&-1&0&0&0&0&0&0\\
28&-1&0&0&0&0&0&0&0&0&0&0&0&0&0&0&b'&0&0&0&0&0&0&0&0\\
29&0&0&1&0&0&0&0&0&0&0&0&0&0&0&0&0&b'&0&0&0&0&0&0&0\\
30&0&0&0&0&0&-1&0&0&0&0&0&0&0&0&0&0&0&b'&0&0&0&0&0&0\\
31&0&0&0&c'&0&0&0&0&0&0&0&0&1&0&0&0&0&0&0&0&0&0&0&0\\
32&0&0&0&0&c'&0&0&0&0&0&0&0&0&1&0&0&0&0&0&0&0&0&0&0\\
33&0&0&0&0&0&c'&0&0&0&0&0&0&0&0&0&1&0&0&0&0&0&0&0&0\\
34&0&0&0&0&0&1&0&0&0&0&0&-1&0&0&0&0&0&0&0&0&0&0&0&0\\
35&0&0&0&0&1&0&0&0&0&1&0&0&0&0&0&0&0&0&0&0&0&0&0&0\\
36&0&0&1&0&0&0&-1&0&0&0&0&0&0&0&0&0&0&0&0&0&0&0&0&0\\
37&0&0&0&0&0&-1&0&0&b&0&0&0&0&0&0&0&0&0&0&0&0&0&0&0\\
38&0&0&0&1&0&0&0&b&0&0&0&0&0&0&0&0&0&0&0&0&0&0&0&0\\
39&-1&0&0&0&0&0&b&0&0&0&0&0&0&0&0&0&0&0&0&0&0&0&0&0\\
40&0&0&c&0&0&0&0&0&0&0&0&1&0&0&0&0&0&0&0&0&0&0&0&0\\
41&0&c&0&0&0&0&0&0&0&0&1&0&0&0&0&0&0&0&0&0&0&0&0&0\\
42&c&0&0&0&0&0&0&0&1&0&0&0&0&0&0&0&0&0&0&0&0&0&0&0\\
43&0&0&0&0&0&0&1&0&0&0&0&0&0&0&0&0&0&0&-1&0&0&0&0&0\\
44&0&0&0&0&0&0&0&1&0&0&0&0&0&0&0&0&0&0&0&0&1&0&0&0\\
45&0&0&0&0&0&0&0&0&0&1&0&0&0&0&0&0&0&0&0&0&0&0&0&-1\\
46&0&0&0&0&0&0&-1&0&0&0&0&0&0&0&0&0&0&0&0&0&0&b'&0&0\\
47&0&0&0&0&0&0&0&0&1&0&0&0&0&0&0&0&0&0&0&0&0&0&b'&0\\
48&0&0&0&0&0&0&0&0&0&0&0&-1&0&0&0&0&0&0&0&0&0&0&0&b'\\
49&0&0&0&0&0&0&0&0&0&c'&0&0&0&0&0&0&0&0&1&0&0&0&0&0\\
50&0&0&0&0&0&0&0&0&0&0&c'&0&0&0&0&0&0&0&0&1&0&0&0&0\\
51&0&0&0&0&0&0&0&0&0&0&0&c'&0&0&0&0&0&0&0&0&0&1&0&0\\
52&0&0&0&0&0&0&0&0&0&0&0&0&0&0&0&0&0&1&0&0&0&0&0&-1\\
53&0&0&0&0&0&0&0&0&0&0&0&0&0&0&0&0&1&0&0&0&0&1&0&0\\
54&0&0&0&0&0&0&0&0&0&0&0&0&0&0&1&0&0&0&-1&0&0&0&0&0\\
55&0&0&0&0&0&0&0&0&0&0&0&0&0&0&0&0&0&-1&0&0&b&0&0&0\\
56&0&0&0&0&0&0&0&0&0&0&0&0&0&0&0&1&0&0&0&b&0&0&0&0\\
57&0&0&0&0&0&0&0&0&0&0&0&0&-1&0&0&0&0&0&b&0&0&0&0&0\\
58&0&0&0&0&0&0&0&0&0&0&0&0&0&0&c&0&0&0&0&0&0&0&0&1\\
59&0&0&0&0&0&0&0&0&0&0&0&0&0&c&0&0&0&0&0&0&0&0&1&0\\
60&0&0&0&0&0&0&0&0&0&0&0&0&c&0&0&0&0&0&0&0&1&0&0&0
\end{array}$$}}

\section{Our specific obstruction class}\label{classsec}
In \cite{CF}, the following theorem is proved.
\begin{thm} Let $X$ be an $n$-dimensional $\Delta$-complex with a single vertex $x_0$, and let $\pi=\pi_1(X,x_0)$.
Let $I\subset\bz[\pi]$ denote the augmentation ideal. The action of $\pi\times \pi$  on $I$ by $(g,h)\cdot\a=g\a h^{-1}$ makes $I$ a $\bz[\pi\times \pi]$-module, defining a local coefficient system $I$ over $X\times X$. If $C_1(X\times X)$ denotes the free abelian group on the set of edges of $X\times X$, then the homomorphism $f:C_1(X\times X)\to I$ defined by
$(e_1,e_2)\mapsto [e_1][e_2]^{-1}-1$ defines an element $\nu\in H^1(X\times X;I)$. Then $\TC(X)< 2n$ iff
$\nu^{2n}=0\in H^{2n}(X\times X;I^{\otimes 2n})$, where $\pi\times\pi$ acts diagonally on $I^{\ot 2n}$.\label{CFthm}
\end{thm}
The discussion of the function $f$ of this theorem in \cite{CF} refers to  \cite[Ch.6:Thm 3.3]{Wh},
and we use the proof of that result for our interpretation of the function.

When $X=K$ is the Klein bottle, we have $\pi=\pi_1(K)=\{a^mb^n\}$ with the multiplication of these elements determined by the relation $ab^{-1}=ba$. Also relevant for us is the element $c=ab^{-1}=ba$.
The ideal  $I$ for us is the free abelian group with basis $\{\a_{m,n}=a^mb^n-1:(m,n)\in\bz\times\bz-\{(0,0)\}\}$.
Using the numbering of the 1-cells $e_i$ of $K\times K$ given in Section \ref{Deltasec}, we obtain that the function $f$
is given as in the following table.

\begin{equation}\begin{array}{c|cccccccc|}\label{ar1}
i&1&2&3&4&5&6&7&8\\
f(e_i)&\a_{1,-1}&\a_{1,0}&\a_{0,1}&\a_{-1,-1}&\a_{-1,0}&\a_{0,-1}&0&\a_{0,1}\\
\hline
\hline
i&9&10&11&12&13&14&15&\\
f(e_i)&\a_{1,-2}&\a_{0,-1}&0&\a_{1,-1}&\a_{-1,-2}&\a_{-1,-1}&0&\\
\hline
\end{array}\end{equation}

\medskip
\noindent For example, $f(e_1)=\a_{1,-1}$ because $c=ab^{-1}$, while $f(e_{12})=\a_{1,-1}$ because the edge from 0 to 2 is $a$, while that from 1 to 2 is $b$.

In \cite[p.500]{EZ}, it is noted that the Alexander-Whitney formula for cup products in simplicial complexes applies also to $\Delta$-complexes.
Our mistake was to overlook the twisting in this formula due to local coefficients. The correct formula, from \cite{St}, is
\begin{equation}\label{twist}(f^p\cup g^{n-p})\la v_0,\ldots,v_n\ra=(-1)^{p(n-p)}f(v_0,\ldots,v_p)\ot \rho_{v_0,v_p}g(v_p,\ldots,v_n),\end{equation}
where $\rho_{i,j}$ is as in (\ref{delta}) and is what we overlooked.

 We apply this to $f^4(\g_j)$, where $f^4=f\cup f\cup f\cup f$ with $f$  the function on 1-cells defined above, and $\g_j$ is any of the 24 4-cells listed in Section \ref{Deltasec}. For example,
\begin{eqnarray}f^4(\g_1)&=&f(00,10)\ot c\,\! f(10,20)\ot a\,\! f(20,21)\ot a\,\! f(21,22)c^{-1}\nonumber\\&=&f(e_1)\ot c\,\! f(e_3)\ot a\,\! f(e_4)\ot a\,\! f(e_6)c^{-1}\nonumber\\
&=&(ab^{-1}-1)\ot c(b-1)\ot a(a^{-1}b^{-1}-1)\ot a(b^{-1}-1)c^{-1}\nonumber\\
&=&(ab^{-1}-1)\ot(a-ab^{-1})\ot(b^{-1}-a)\ot(1-b^{-1}).\label{f4}\end{eqnarray}
  Formula (\ref{f4}) equals (T19) in the expansion of the obstruction class in Section 3.2 of \cite{CV}, after adjusting for different notation. They use $x,y$ for our $a,b$, and they write their classes as $y^nx^m$ rather than our $a^mb^n$. The relations in $\pi_1(K)$ must be used to compare these.

For all of our 4-cells, consecutive vertices are constant in one factor, and so only $f(e_i)$ for $i\le 6$ are relevant for $f^4$.
Once we knew about incorporating the twisting in (\ref{twist}), our obstruction class exactly agreed with the class in Section 3.2 of \cite{CV}. Prior to this realization, we worked with a different obstruction class for several months.

\section{Approaches to proving that our class is nonzero}

As just noted, once (\ref{twist}) was understood, our obstruction class agreed with the obstruction class of \cite{CV}. They successfully showed that this was nonzero in $H^4(K\times K;I^{\ot 4})$.
From our viewpoint, the relations were of the form
 \begin{eqnarray*}&(1)&b\a_{m_1,n_1}\ot b\a_{m_2,n_2}\ot b\a_{m_3,n_3}\ot b\a_{m_4,n_4}-\a_{m_1,n_1}\ot\a_{m_2,n_2}\ot \a_{m_3,n_3}\ot \a_{m_4,n_4},\\
 &(2)&\a_{m_1,n_1}b^{-1}\ot \a_{m_2,n_2}b^{-1}\ot \a_{m_3,n_3}b^{-1}\ot \a_{m_4,n_4}b^{-1}\\
 &&\quad-\a_{m_1,n_1}\ot\a_{m_2,n_2}\ot \a_{m_3,n_3}\ot \a_{m_4,n_4},\\
 &(3)&c\a_{m_1,n_1}\ot c\a_{m_2,n_2}\ot c\a_{m_3,n_3}\ot c\a_{m_4,n_4}+\a_{m_1,n_1}\ot\a_{m_2,n_2}\ot \a_{m_3,n_3}\ot \a_{m_4,n_4},\\
 &(4)&\a_{m_1,n_1}c^{-1}\ot \a_{m_2,n_2}c^{-1}\ot \a_{m_3,n_3}c^{-1}\ot \a_{m_4,n_4}c^{-1}\\
 &&\quad+\a_{m_1,n_1}\ot\a_{m_2,n_2}\ot \a_{m_3,n_3}\ot \a_{m_4,n_4},\end{eqnarray*}
 and we hoped to show our class could not be reduced to 0 mod these relations.

The factors $b\a_{m,n}$, $\a_{m,n}b^{-1}$, $c\a_{m,n}$, and $\a_{m,n}c^{-1}$ appearing above are, respectively,
 \begin{eqnarray*}1&:&\a_{m,n+(-1)^m}-\a_{0,1},\\
 2&:&\a_{m,n-1}-\a_{0,-1},\\
 3&:&\a_{m+1,n+(-1)^{m+1}}-\a_{1,-1},\\
 4&:&\a_{m-1,-n-1}-\a_{-1,-1}.\end{eqnarray*}
 The relations are much simpler in $\Z[\pi]^{\ot4}$, and we had a nice form for those. For example,
 $$b a^{m_1}b^{n_1}\ot \cdots\ot b a^{m_4}b^{n_4}\sim a^{m_1}b^{n_1}\ot \cdots\ot a^{m_4}b^{n_4}$$
 is much simpler than (1). Prior to our understanding of the twisting in (\ref{twist}), we had been able to show our class was not 0 mod the $\Z[\pi]^{\ot4}$-relations. Once we were made aware of the twisting,
 we showed our corrected class was 0 mod these relations. Then it was pointed out to us (by Cohen and Vandembroucq)
that it followed from the work in \cite{CF} that the class would necessarily be zero with coefficients in $\Z[\pi]^{\ot4}$. Our methods lent no insight toward showing it nonzero with coefficients $I^{\ot4}$.

We wish to thank the authors of \cite{CV} for their help in our understanding of this project.
Still, we feel that our approach to obtaining the obstruction class, quite different from theirs, is worth publicizing.

 \def\line{\rule{.6in}{.6pt}}

\end{document}